# ON THE UNIQUENESS PROBLEMS OF ENTIRE FUNCTIONS AND THEIR LINEAR DIFFERENTIAL POLYNOMIALS*

QI HAN AND HONG-XUN YI

### Abstract

The uniqueness problems on transcendental meromorphic or entire functions sharing at least two values with their derivatives or linear differential polynomials have been studied and many results have been obtained. In this paper, we study a transcendental entire function $f(z)$ that shares a non-zero polynomial $a(z)$ with $f'(z)$, together with its linear differential polynomials of the form: $L[f] = a_2(z)f''(z) + a_3(z)f'''(z) + \cdots + a_m(z)f^{(m)}(z)$ $(a_m(z) \not\equiv 0)$, where the coefficients $a_k(z)$ $(k = 2, 3, \ldots, m)$ are rational functions.

## 1. Introduction and main results

In this paper, a meromorphic function will always mean meromorphic in the complex plane $\mathbf{C}$. We adopt the standard notations in the Nevanlinna theory of value distribution of meromorphic functions such as *the characteristic function* $T(r, f)$, *the proximity function* $m(r, f)$ and *the counting function* $N(r, f)$ ($\overline{N}(r, f)$) of poles. For any non-constant meromorphic function $f$, we denote by $S(r, f)$ any quantity satisfying $S(r, f) = o(T(r, f))$, possibly outside a set of finite linear measure in $\mathbf{R}$. Furthermore, the exceptional set is not necessarily the same at each occurrence.

Let $f$ be a non-constant meromorphic function on $\mathbf{C}$, and let $a \in \mathcal{M}_f(\mathbf{C})$ be any meromorphic function whose characteristic function satisfies $T(r, a) = S(r, f)$. Then, we call $a$ a small function of $f$ and $\mathcal{M}_f(\mathbf{C})$ the set of small functions of $f$ over $\mathbf{C}$. Obviously, $\mathcal{M}_f(\mathbf{C})$ is a field and contains $\mathbf{C}$.

Let $f$ be a non-constant meromorphic function on $\mathbf{C}$, let $a \in \mathcal{M}_f(\mathbf{C})$ be a small function of $f$, and let $\mu_f^a(z)$ denote the multiplicity of $f - a = 0$ at $z$. We define the *a*-valued *divisor* of $f$ by

$$E_f(a) = \{(\mu_f^a(z), z) \in \mathbf{N} \times \mathbf{C} \,|\, z \in \mathbf{C}\},$$

2000 *Mathematics Subject Classification*: Primary 30D35. Secondly 30D20, 34M05.
*Supported by the NNSF of P. R. China and the NSF of Shandong Province.
Received July 1, 2005; revised October 2, 2006.

and the *preimage* of $a$ under $f$ by

$$\overline{E}_f(a) = f^{-1}(a) = \{z \in \mathbf{C} \mid \mu_f^a(z) > 0\}.$$

Let $g$ be another non-constant meromorphic function on $\mathbf{C}$. We say that $f$ and $g$ share $a$ CM (resp., IM) whenever $E_f(a) = E_g(a)$ (resp., $\overline{E}_f(a) = \overline{E}_g(a)$), which means counting (resp., ignoring) multiplicities. Also, we say that $f$ and $g$ share the value $\infty$ CM (resp., IM) whenever $1/f$ and $1/g$ share the value $0$ CM (resp., IM).

It is well-known that a pair of non-constant meromorphic functions $f$ and $g$ would be identically equal to each other if $f$ and $g$ share five distinct values IM. That is the famous *Nevanlinna's Five-value Theorem*. Also, the *Nevanlinna's Four-value Theorem* states that if a pair of non-constant meromorphic functions $f$ and $g$ share four distinct values CM, then $f$ is a bilinear transformation of $g$. The condition that "$f$ and $g$ share four distinct values CM" has been weakened to "$f$ and $g$ share two distinct values CM and other two distinct values IM" by G. Gundersen, i.e., "2CM + 2IM = 4CM" holds.

We refer the reader to [10] and [21] for those foregoing notations and results.

If we focus on the value-sharing problems of $f$, a non-constant meromorphic function, and its derivatives or its linear differential polynomials, then the number of the shared values could be reduced. For example, if $f$ shares two finite, distinct values CM with $f^{(k)}$ ($k \in \mathbf{N}$), then $f \equiv f^{(k)}$ (see [7], [9] and [19]); if $f$ shares three finite, distinct values IM with $P[f] := b_0 f + b_1 f' + \cdots + b_n f^{(n)}$ ($b_n \neq 0$), then $f \equiv P[f]$ (see [8], [15] and [20]), where $b_l$ are constants for $l = 0, 1, \ldots, n$. If $f$ is entire and shares two finite, distinct values with its linear differential polynomials of similar forms to $P[f]$ but with small function coefficients, say, $b_l \in \mathscr{M}_f(\mathbf{C})$ ($l = 0, 1, \ldots, n$), many results have been obtained (see [1], [14] and [16]).

It were Jank-Mues-Volkmann who first studied the problem of sharing only one finite non-zero value between an entire function $f$ and its derivatives. As a matter of fact, in 1986, they proved the following theorem in [12].

THEOREM A. *Let $f$ be a non-constant entire function. If $f$ and $f'$ share a finite, non-zero value $a \in \mathbf{C}$ IM, and $f'' = a$ whenever $f = a$, then $f \equiv f'$.*

In fact, from the hypothesis of Theorem A, it could be easily seen that the value $a$ is shared by $f$ and $f'$ CM.

In 2001, Li-Yang obtained the following result in [17].

THEOREM B. *Let $f$ be a non-constant entire function, let $a \in \mathbf{C}$ be a finite, non-zero value, and let $k(\geq 2)$ be a positive integer. If $f$, $f'$ and $f^{(k)}$ share the value $a$ CM, then $f$ assumes the form*

$$f(z) = ce^{bz} + \frac{a(b-1)}{b},$$

*where $b$, $c$ are two non-zero constants and $b$ satisfies $b^{k-1} = 1$.*

In [2] and [3], Chang-Fang generalized Theorems A and B. In effect, their results state that if we replace the value $a \in \mathbf{C}\setminus\{0\}$ in Theorem A with any element $a \in \mathscr{M}_f(\mathbf{C})$ satisfying $a \not\equiv a'$ and the value $a \in \mathbf{C}\setminus\{0\}$ in Theorem B with any element $a \in \mathscr{M}_f(\mathbf{C})\setminus\mathbf{C}$; and if we retain all the other assumptions in Theorems A and B, respectively, then $f \equiv f'$. However, the methods of their proofs were similar to those original ones in [12] and [17], respectively.

In this paper, we investigate the uniqueness problem of a transcendental entire function $f$, its first derivative $f'$ and its linear differential polynomials of a certain form. In fact, by employing a different method from those of Theorems A or B but more or less inspired by that of Theorem B, we obtain the following result.

THEOREM 1. *Let $f$ be a transcendental entire function, and let*

$$L[f] := a_2(z)f''(z) + a_3(z)f'''(z) + \cdots + a_m(z)f^{(m)}(z) \quad (a_m(z) \not\equiv 0)$$

*be a linear differential polynomial in $f$ with rational coefficients $a_k(z)$ ($k = 2, 3, \ldots, m$). If $f$ and $f'$ share a non-zero polynomial $a$ of degree $p$ CM, and $L[f] = a$ whenever $f = a$, then $f$ assumes the form*

$$f(z) = ce^{bz} + (b-1)\Phi(z),$$

*where $b$, $c$ are two non-zero constants and $\Phi(z)$ is defined as*

$$\Phi(z) := \sum_{j=1}^{p+1} \frac{a^{(j-1)}(z)}{b^j}.$$

*In particular, $\Phi(z)$ satisfies the equation*

$$\Psi(z) \equiv \frac{a(z) - (b-1)L[\Phi]}{a(z) - (b-1)\Phi(z)},$$

*where $\Psi(z)$ is defined as $\Psi(z) := \sum_{k=2}^{m} a_k(z)b^k$.*

*Remark* 1. Obviously, Theorem 1 generalizes Theorems A and B, since $f^{(k)}$ ($k \in \mathbf{N}, k \geq 2$) is a special form of $L[f]$ and any non-zero constant is a polynomial of degree 0. Also, the assumptions that "$f$, $f'$ share $a$ CM" and "$L[f] = a$ whenever $f = a$" are weaker than the assumption that "$f$, $f'$ and $L[f]$ share $a$ CM". Furthermore, it does not seem feasible to apply the methods in [12] or [17] directly to the assumptions of Theorem 1 to get its conclusion.

*Remark* 2. Let $f$ be a non-constant meromorphic function on $\mathbf{C}$, and let $n_0, n_1, \ldots, n_k$ be $k+1$ non-negative integers. We call

$$M[f] := f^{n_0}(f')^{n_1} \cdots (f^{(k)})^{n_k}$$

a *monomial* in $f$ of *degree* $\gamma_M := n_0 + n_1 + \cdots + n_k$. Now, let $M_1[f], M_2[f], \ldots, M_l[f]$ be $l$ monomials in $f$ of degree $\gamma_{M_j}$ ($j = 1, 2, \ldots, l$), respectively; and let $b_1, b_2, \ldots, b_l$ be $l$ elements in $\mathscr{M}_f(\mathbf{C})$. We call

$$M^*[f] := b_1 M_1[f] + b_2 M_2[f] + \cdots + b_l M_l[f]$$

a *differential polynomial* in $f$ of *total degree* $\gamma_{M^*} := \max\{\gamma_{M_1}, \gamma_{M_2}, \ldots, \gamma_{M_l}\}$.

In particular, $M^*[f]$ ia called *linear* if $\gamma_{M^*} = 1$ (see [6], [10] or [21]).

## 2. Some lemmas

LEMMA 1 (see [5] and [6]). *Let $f$ be a non-constant meromorphic function, and let $Q^*[f]$ and $Q[f]$ be two differential polynomials in $f$ with arbitrary meromorphic coefficients $q_1^*, q_2^*, \ldots, q_s^*$ and $q_1, q_2, \ldots, q_t$, respectively. Furthermore, let $P(f)$ be a non-constant polynomial of $f$ of degree $p$ whose coefficients are small functions of $f$. If the total degree $\gamma_Q$ of $Q$ is at most $p$, then from the following equation*

$$P(f)Q^*[f] \equiv Q[f],$$

*we have*

$$m(r, Q^*[f]) \leq \sum_{u=1}^{s} m(r, q_u^*) + \sum_{v=1}^{t} m(r, q_v) + S(r, f).$$

LEMMA 2 (see [13, p.p. 58, Remark 1]). *Let $f$ be a solution of the following homogeneous differential equation*

$$a_n(z) f^{(n)}(z) + a_{n-1}(z) f^{(n-1)}(z) + \cdots + a_1(z) f'(z) + a_0(z) f(z) = 0,$$

*where the coefficients $a_0(z), \ldots, a_n(z)$ are polynomials and are not all identically equal to zero. Then, $f$ is an entire function of finite order.*

The order of a meromorphic function, say, $f$, is defined as $\sigma(f) := \limsup_{r \to +\infty} \frac{\log T(r, f)}{\log r}$. Then, the well-known *Chuang's inequality* (see [4]) and the *lemma of logarithmic derivative* (see [10, p.p. 34–42]) yield $\sigma(f) = \sigma(f^{(k)})$ for all positive integers $k \in \mathbf{N}$.

LEMMA 3. *Let $P(z)$ and $Q(z)$ be two polynomials with $P(z) \not\equiv 0$ and $Q(z) \in \mathbf{C}[z] \backslash \mathbf{C}$. Then, every entire solution $F(z)$ of the following differential equation*

$$F'(z) - e^{Q(z)} F(z) = P(z)$$

*has infinite order.*

The proof of Lemma 3 will be given in the next part since it is relatively long and of interest in itself.

LEMMA 4 (see [18] or [21]). *Let $f$ be a non-constant meromorphic function, and let $R(f) = \sum_{s=0}^{m} a_s f^s / \sum_{t=0}^{n} b_t f^t$ be an irreducible rational function in $f$ with small function coefficients $\{a_s\}_{s=0}^{m}$ and $\{b_t\}_{t=0}^{n}$ $(a_m \cdot b_n \not\equiv 0)$. Then,*

$$T(r, R(f)) = \max\{m, n\} T(r, f) + S(r, f).$$

## 3. Proof of Lemma 3

On the contrary, we assume that there exists some entire solution $F(z)$ of the following differential equation

$$(3.1) \qquad F'(z) - e^{Q(z)} F(z) = P(z)$$

with finite order, say, $\sigma \geq 0$. The classical Wiman-Valiron estimate (see [11], also see [13] for an analogous form) states that, for the transcendental entire function $F(z)$, and for any positive real number $\varepsilon$ with $0 < \varepsilon < 1/8$, if we let $z$ with $|z| = r$ be such that the following inequality

$$(3.2) \qquad |F(z)| > M(r, F)(v(r, F))^{-1/8+\varepsilon}$$

holds, then there exists a set $\mathbf{E} \subset \mathbf{R}(1, \infty)$ with finite logarithmic measure such that the following estimate

$$(3.3) \qquad \frac{F^{(k)}(z)}{F(z)} = \left(\frac{v(r, F)}{z}\right)^k (1 + O(v(r, F)^{-1/8+\varepsilon}))$$

holds for all $k \in \mathbf{N}$ and all $r \notin \mathbf{E} \cup \mathbf{R}[0, 1]$, where $M(r, F)$ denotes the *maximal module* of $F(z)$ at $|z| = r$. As to the term $v(r, F)$, called the *central index* of $F(z)$, is defined as $\max_{n \geq 0}\{n \mid |a_n| r^n = \mu(r, F)\}$ with the *maximum term* of $F(z)$ defined as $\mu(r, F) := \max_{n \geq 0}\{|a_n| r^n\}$ if the Taylor expansion of $F(z)$ is $F(z) = \sum_{n=0}^{\infty} a_n z^n$ with center 0 and radius $r$. Furthermore, we have $\sigma = \limsup_{r \to +\infty} \frac{\log^+ v(r, F)}{\log r}$ with $\log^+ x := \max\{\log x, 0\}$ for $x \geq 0$.

Suppose that the degree of $P(z)$ is $l$. Then, differentiating both sides of Eq. (3.1) $l + 1$ times yields the following homogeneous differential equation

$$(3.4) \qquad F^{(l+2)}(z) - e^{Q(z)}(F^{(l+1)}(z) + Q_l(z) F^{(l)}(z) + \cdots + Q_0(z) F(z)) = 0,$$

where the coefficients $Q_j(z)$ $(j = 0, 1, \ldots, l)$ are polynomials written in terms of $Q(z)$, its derivatives and their combinations.

Let $\{z_n\}_{n=1}^{\infty}$ be a sequence of complex numbers with $|z_n| \to \infty$ as $n \to \infty$ such that it satisfies Eq. (3.2) and Eq. (3.3) simultaneously, and let $Q(z)$ be of the form $Q(z) = q_m z^m + q_{m-1} z^{m-1} + \cdots + q_0$ with $q_m \neq 0$. Now, suppose that there exists a subsequence $\{z_{n_k}\}_{k=1}^{\infty}$ of $\{z_n\}_{n=1}^{\infty}$ such that

$$(3.5) \qquad |e^{Q(z_{n_k})}| \sim \exp\{\Re(q_m z_{n_k}^m)(1 + o(1))\}$$

holds, and such that $\Re(q_m z_{n_k}^m) > 0$ holds for sufficiently large $k$'s. Substituting Eq. (3.3) and Eq. (3.5) into Eq. (3.4) written in the following form

(3.6) $$\frac{F^{(l+2)}(z)}{F(z)} - e^{Q(z)}\left(\frac{F^{(l+1)}(z)}{F(z)} + Q_l(z)\frac{F^{(l)}(z)}{F(z)} + \cdots + Q_0(z)\right) = 0$$

yields a contradiction since we assume that $F(z)$ has finite order.

On the other hand, if there exists no such subsequence $\{z_{n_k}\}_{k=1}^{\infty}$ of $\{z_n\}_{n=1}^{\infty}$ for which Eq. (3.5) holds could be found with $\Re(q_m z_{n_k}^m) > 0$, then we should have either $\Re(q_m z_{n_k}^m) < 0$ or $\Re(q_m z_{n_k}^m) = 0$ for an infinite number of points from $\{z_n\}_{n=1}^{\infty}$. If $\Re(q_m z_{n_k}^m) < 0$ holds for an infinite number of $k$'s with $k \to \infty$, then substituting Eq. (3.3) and Eq. (3.5) into Eq. (3.6) yields a contradiction again.

Finally, if $\Re(q_m z_{n_k}^m) = 0$ holds for all subsequences $\{z_{n_k}\}_{k=1}^{\infty}$ of $\{z_n\}_{n=1}^{\infty}$, then we would consider the next term $\Re(q_{m-1} z_{n_k}^{m-1})$ instead. If $\Re(q_{m-1} z_{n_k}^{m-1}) = 0$ holds for all subsequences $\{z_{n_k}\}_{k=1}^{\infty}$ of $\{z_n\}_{n=1}^{\infty}$ again, then we would repeat the argument to yet another lower term of $Q(z)$ and so on and so forth. We could eventually derive a contradiction by repeating the argument, which completes the proof. □

## 4. Proof of Theorem 1

We now define the function $\alpha$ to be

(4.1) $$\alpha := \frac{(a-a')(L[f] - L[a]) - (a - L[a])(f' - a')}{f - a}.$$

By the *lemma of logarithmic derivative*, we have $m(r, \alpha) = S(r, f)$. It is not difficult to see that the poles of $\alpha$ arise from the zeros of $a - a'$, since $f$ and $f'$ share $a$ CM, and the poles of $a_k$ ($k = 2, 3, \ldots, m$). Since $f$ is transcendental, we get

(4.2) $$T(r, \alpha) = m(r, \alpha) + O(\log r) = S(r, f).$$

We now distinguish the following two cases.

Case (i). $\alpha$ is rational.

Write $\alpha$ in the form $\alpha(z) = p(z)/q(z)$, where $p(z)$ and $q(z)$ are two polynomials with no common factors. It is not difficult to see that $f$ satisfies a linear differential equation with polynomial coefficients by Eq. (4.1), and hence has finite order by the conclusion of Lemma 2 and the fact that $\sigma(f) = \sigma(f^{(k)})$ ($k \in \mathbf{N}$).

Since $f$ and $f'$ share $a$ CM, we may suppose that

(4.3) $$\frac{f'(z) - a(z)}{f(z) - a(z)} = e^{Q(z)},$$

where $Q(z)$ is an entire function. Rewrite Eq. (4.3) as

$$(f'(z) - a'(z)) - (f(z) - a(z))e^{Q(z)} = a(z) - a'(z).$$

We now apply the conclusion of Lemma 3 to the above equation to derive that $Q(z)$ is a constant. Thus, there exists a constant $b \in \mathbf{C}\setminus\{0\}$ such that

$$\frac{f'(z) - a(z)}{f(z) - a(z)} = b.$$

Solving the above non-homogeneous linear differential equation yields

(4.4) $$f(z) = ce^{bz} + (b-1)\Phi(z),$$

where $c(\neq 0)$ is a constant and $\Phi(z)$ is defined as $\Phi(z) := \sum_{j=1}^{p+1} \frac{a^{(j-1)}(z)}{b^j}$. So,

(4.5) $$T(r, f) = T(r, e^{bz}) + S(r, f),$$
$$f^{(k)}(z) = cb^k e^{bz} + (b-1)\Phi^{(k)}(z) \quad (k \geq 1),$$

and

(4.6) $$L[f(z)] = ce^{bz}\Psi(z) + (b-1)L[\Phi(z)],$$

where $\Psi(z)$ is defined as $\Psi(z) := \sum_{k=2}^{m} a_k(z) b^k$.

Since $e^{bz}$ has two Picard values 0 and $\infty$, noting Eq. (4.4) and Eq. (4.5), we see that $f$ has two Picard "small functions" $(b-1)\Phi$ and $\infty$, i.e.,

$$N(r, f) = S(r, f), \quad N\left(r, \frac{1}{f - (b-1)\Phi}\right) = S(r, f).$$

Obviously, $(b-1)\Phi(z) \not\equiv a(z)$; otherwise, we would have $-1 = 0$, a contradiction. We now define the function $f^*$ to be $f^* := \frac{f - a}{f - (b-1)\Phi}$. Hence,

(4.7) $$T(r, f) = T(r, f^*) + S(r, f),$$

and

(4.8) $$N\left(r, \frac{1}{f - a}\right) = N\left(r, \frac{1}{f^*}\right) + S(r, f).$$

Furthermore, $N(r, f^*) = N\left(r, \frac{1}{f - (b-1)\Phi}\right) = S(r, f)$ and $N\left(r, \frac{1}{f^* - 1}\right) = N(r, f) + O(\log r) = S(r, f)$. By Eq. (4.7), we also have $S(r, f) = S(r, f^*)$.

We now apply the *Second Main Theorem* to $f^*$ with 0, 1 and $\infty$ to conclude

$$T(r, f^*) = N\left(r, \frac{1}{f^*}\right) + S(r, f^*).$$

Hence, combining the above equation with Eq. (4.7) and Eq. (4.8) yields

(4.9) $$T(r, f) = N\left(r, \frac{1}{f - a}\right) + S(r, f).$$

From Eq. (4.9), we see that $f - a$ has infinitely many zeros. Take $z_a$ to be a zero of $f - a$. Since $L[f] = a$ whenever $f = a$, from Eq. (4.4) and Eq. (4.6), we derive

$$a(z_a) = ce^{bz_a} + (b-1)\Phi(z_a),$$
$$a(z_a) = ce^{bz_a}\Psi(z_a) + (b-1)L[\Phi(z_a)].$$

Eliminating $ce^{bz_a}$ from the above two equations yields

$$a(z_a) - (a(z_a) - (b-1)\Phi(z_a))\Psi(z_a) - (b-1)L[\Phi(z_a)] = 0.$$

If $a(z) - (a(z) - (b-1)\Phi(z))\Psi(z) - (b-1)L[\Phi]$ is not identically equal to zero, then we derive,

$$N\left(r, \frac{1}{f-a}\right) \le N\left(r, \frac{1}{a - (a-(b-1)\Phi)\Psi - (b-1)L[\Phi]}\right)$$
$$\le T\left(r, \frac{1}{a - (a-(b-1)\Phi)\Psi - (b-1)L[\Phi]}\right) = O(\log r) = S(r, f),$$

which contradicts Eq. (4.9). Hence,

(4.10) $$a(z) \equiv (a(z) - (b-1)\Phi(z))\Psi(z) + (b-1)L[\Phi],$$

which can be rewritten as $\Psi(z) \equiv \dfrac{a(z) - (b-1)L[\Phi]}{a(z) - (b-1)\Phi(z)}$.

CASE (ii). $\alpha$ is transcendental.
Rewrite Eq. (4.1) as

$$f - a = \frac{1}{\alpha}((a-a')(L[f] - L[a]) - (a - L[a])(f' - a')).$$

Differentiating the above equation yields

(4.11) $$f' - a' = \left(\frac{1}{\alpha}\right)'((a-a')(L[f] - L[a]) - (a - L[a])(f' - a'))$$
$$+ \frac{1}{\alpha}((a'-a'')(L[f] - L[a]) + (a-a')(L'[f] - L'[a])$$
$$- (a' - L'[a])(f' - a') - (a - L[a])(f'' - a'')).$$

A routine calculation leads to

(4.12) $$\left(1 + (a - L[a])\left(\frac{1}{\alpha}\right)' + (a' - L'[a])\frac{1}{\alpha}\right)(f' - a)$$
$$= \sum_{k=2}^{m+1} a_k^*(f^{(k)} - a^{(k-1)})$$
$$+ \left\{(a' - a) + \left(\frac{1}{\alpha}\right)'\left((a-a')\sum_{k=2}^{m} a_k a^{(k-1)} - a(a-a')\right)\right.$$

$$+\frac{1}{\alpha}\left((a-a')\left(\sum_{k=2}^m a'_k a^{(k-1)} - a' + \sum_{k=2}^m a_k a^{(k)}\right)\right.$$
$$\left.+(a'-a'')\left(\sum_{k=2}^m a_k a^{(k-1)} - a\right)\right)\bigg\},$$

where $a_k^*$ $(k=2,3,\ldots,m+1)$ are linear combinations of $a_j$ $(j=2,3,\ldots,m)$, $1/\alpha$ and their derivatives, $a$ and $a-a'$ etc..

Since $a \not\equiv a'$, we claim that the following expression

(4.13) $$(a'-a) + \left(\frac{1}{\alpha}\right)'\left((a-a')\sum_{k=2}^m a_k a^{(k-1)} - a(a-a')\right)$$
$$+\frac{1}{\alpha}\left((a-a')\left(\sum_{k=2}^m a'_k a^{(k-1)} - a' + \sum_{k=2}^m a_k a^{(k)}\right)\right.$$
$$\left.+(a'-a'')\left(\sum_{k=2}^m a_k a^{(k-1)} - a\right)\right)$$

could not be identically equal to zero. If not, we rewrite Eq. (4.13) as

$$a - a' \equiv \frac{1}{\alpha}\left((a-a')\left(\sum_{k=2}^m a'_k a^{(k-1)} - a' + \sum_{k=2}^m a_k a^{(k)}\right)\right.$$
$$\left.+(a'-a'')\left(\sum_{k=2}^m a_k a^{(k-1)} - a\right)\right)$$
$$+\left(\frac{1}{\alpha}\right)'\left((a-a')\sum_{k=2}^m a_k a^{(k-1)} - a(a-a')\right).$$

Dividing $(a-a')$ on both sides of the above equation yields

(4.13*) $$1 \equiv \frac{1}{\alpha}\left(\left(\sum_{k=2}^m a'_k a^{(k-1)} - a' + \sum_{k=2}^m a_k a^{(k)}\right) + \left(\frac{a'-a''}{a-a'}\right)\left(\sum_{k=2}^m a_k a^{(k-1)} - a\right)\right)$$
$$+\left(\frac{1}{\alpha}\right)'\left(\sum_{k=2}^m a_k a^{(k-1)} - a\right).$$

Multiplying $\alpha$ on both sides of Eq. (4.13*), plus the fact that $N(r,\alpha) = O(\log r)$ from Eq. (4.2) and the *lemma of logarithmic derivative*, yields

$$T(r,\alpha) = m(r,\alpha) + O(\log r) \leq S(r,\alpha) + O(\log r),$$

which contradicts the assumption that $\alpha$ is transcendental.

Let $\beta$ denote the expression of Eq. (4.13) which is not identically equal to zero. Dividing both sides of Eq. (4.12) by $\beta$ yields

$$(4.12^*) \quad 1 = \frac{\left(1 + (a - L[a])\left(\frac{1}{\alpha}\right)' + (a' - L'[a])\frac{1}{\alpha}\right)(f' - a)}{\beta}$$

$$- \frac{\sum_{k=2}^{m+1} a_k^*(f^{(k)} - a^{(k-1)})}{\beta}.$$

Dividing $(f' - a)$ on both sides of Eq. (4.12*), plus the *lemma of logarithmic derivative* and the fact that $a$, $a_k$ $(k = 2, 3, \ldots, m)$, $\alpha$ and thus $\beta$ are all small functions of $f$, yields

$$(4.14) \quad m\left(r, \frac{1}{f' - a}\right) = S(r, f).$$

Since $f$ is entire, we have

$$T(r, f) = m(r, f) + O(1)$$
$$= m\left(r, a + \frac{1}{\alpha}((a - a')(L[f] - L[a]) - (a - L[a])(f' - a'))\right) + O(1)$$
$$\leq m(r, f') + S(r, f) \leq T(r, f') + S(r, f) \leq T(r, f) + S(r, f);$$

that is,

$$(4.15) \quad T(r, f) = T(r, f') + S(r, f).$$

By Eq. (4.14), Eq. (4.15) and the assumption that $f$ and $f'$ share $a$ CM, we get

$$(4.16) \quad m\left(r, \frac{1}{f - a}\right) = T(r, f) - N\left(r, \frac{1}{f - a}\right) + S(r, f)$$
$$= T(r, f') - N\left(r, \frac{1}{f' - a}\right) + S(r, f)$$
$$= m\left(r, \frac{1}{f' - a}\right) + S(r, f) = S(r, f).$$

Thus, Eq. (4.9) still holds.

We now define $\psi$ to be

$$(4.17) \quad \psi := \frac{f' - a}{f - a}.$$

Obviously, by the *lemma of logarithmic derivative* and Eq. (4.16), plus the assumption that $f$ and $f'$ share $a$ CM, we see that $\psi$ is a small entire function of $f$, i.e., $T(r, \psi) = S(r, f)$. Rewrite Eq. (4.17) as

$$(4.18) \quad f' = \psi f + a(1 - \psi) = \lambda_1 f + \mu_1,$$

where $\lambda_1$ and $\mu_1$ are defined as $\lambda_1 := \psi$, $\mu_1 := a(1 - \psi)$.

We repeat the above argument $(k-1)$-times by differentiating Eq. (4.18) repeatedly and by substituting the differential equations into earlier ones to conclude

(4.19) $$f^{(k)} = \lambda_k f + \mu_k \quad (k = 1, 2, \ldots),$$

where $\lambda_k$ and $\mu_k$ are small entire functions with the following recurrence formulas:

$$\lambda_k = \lambda'_{k-1} + \lambda_1 \lambda_{k-1}, \quad \mu_k = \mu'_{k-1} + \mu_1 \lambda_{k-1} \quad (k = 2, 3, \ldots).$$

It is not difficult, by applying induction in the number $k$ to $\lambda_k$ and $\mu_k$, plus the knowledge that differentiation never increases the total degree of a differential polynomial, to deduce that

(4.20) $$\lambda_k = \psi^k + P_{k-1}[\psi], \quad \mu_k = -a\psi^k + P^*_{k-1}[\psi] \quad (k = 1, 2, \ldots),$$

(4.21) $$\mu_{k+1} + a\lambda_{k+1} = (a - a')\psi^k + Q_{k-1}[\psi] \quad (k = 1, 2, \ldots),$$

where $P_{k-1}[\psi]$, $P^*_{k-1}[\psi]$ and $Q_{k-1}[\psi]$ are differential polynomials in $\psi$ with polynomial coefficients and $\max\{\gamma_{P_{k-1}}, \gamma_{P^*_{k-1}}, \gamma_{Q_{k-1}}\} \leq k - 1$. Here, $\gamma_{P_{k-1}}$ denotes the total degree of $P_{k-1}[\psi]$ and $\gamma_{P^*_{k-1}}$, $\gamma_{Q_{k-1}}$ are similarly defined.

Since $L[f] = a$ whenever $f = a$, by Eq. (4.9), Eq. (4.19), the expression of $L[f]$ and the fact that $\psi$, $a$, $a_k$, $\lambda_k$ and $\mu_k$ $(k = 2, 3, \ldots, m)$ are all small functions of $f$, plus analogous reasoning to that of Eq. (4.10), we have

(4.22) $$a \equiv \sum_{k=2}^{m} a_k(a\lambda_k + \mu_k) \equiv \sum_{k=2}^{m} a_k((a - a')\psi^{k-1} + Q_{k-2}[\psi])$$
$$\equiv (a - a')a_m \psi^{m-1} + Q^*_{m-2}[\psi].$$

Obviously, $Q^*_{m-2}[\psi] \not\equiv a$; otherwise, we would have $\psi \equiv 0$, a contradiction. We now apply the conclusion of Lemma 1 to Eq. (4.22) to conclude

$$T(r, \psi) = m(r, \psi) + O(1) \leq S(r, \psi) + O(\log r),$$

which implies that $\psi$ is a polynomial without zeros, and hence a non-zero constant, say, $b \in \mathbf{C}/\{0\}$. Thus, Eq. (4.4), Eq. (4.5) and Eq. (4.6) still hold.

On the other hand, by Eq. (4.4) and Eq. (4.6), Eq. (4.1) can be written as

(4.23) $$\alpha = \frac{ce^{bz}((a - a')\Psi - b(a - L[a]))}{ce^{bz} + (b - 1)\Phi - a}$$
$$+ \frac{(a - a')((b - 1)L[\Phi] - L[a]) - (a - L[a])((b - 1)\Phi' - a')}{ce^{bz} + (b - 1)\Phi - a}.$$

We now apply the conclusion of Lemma 4 to Eq. (4.23), noting Eq. (4.5), to obtain that either $T(r, \alpha) = O(\log r)$ or $T(r, \alpha) = T(r, f) + O(\log r)$ dependent upon whether $\alpha$ is reducible or not. However, it contradicts the assumption that $\alpha$ is a transcendental small function of $f$.

This contradiction finishes the proof of Theorem 1. □

*Acknowledgement.* The first author would like to express his genuine gratitude to his parents for their love and support. Both the authors are indebted to the referee for the careful reading through of the paper with valuable comments and suggestions, especially for the essential advices on great simplification of the proof of Lemma 3.

Qi Han
DEPARTMENT OF MATHEMATICS
SHANDONG UNIVERSITY
250100 JINAN, SHANDONG
P. R. CHINA
E-mail: k.l.han@mail.sdu.edu.cn

CURRENT ADDRESS:
DEPARTMENT OF MATHEMATICS
YAMAGATA UNIVERSITY
1-4-12 KOJIRAKAWA-MACHI, YAMAGATA-SHI 990-8560
JAPAN
E-mail: han@sci.kj.yamagata-u.ac.jp

Hong-Xun Yi
DEPARTMENT OF MATHEMATICS
SHANDONG UNIVERSITY
JINAN, SHANDONG 250100
PEOPLE'S REPBLIC OF CHINA
E-mail: hxyi@sdu.edu.cn